\documentclass[12pt]{article}
\usepackage{amssymb}
\usepackage{latexsym}
\usepackage{amsmath,amscd}
\usepackage{slashed}

\def\part#1{\frac{\partial\phantom{q}}{\partial#1}}

\newenvironment{rmk}{\begin{trivlist}\item[]{\bf Remark:} }
{\end{trivlist}}

\newenvironment{prf}{\begin{trivlist}\item[]{\bf Proof:} }
{\hfill $\Box$ \end{trivlist}} 
\newtheorem{thm}{Theorem}

\newtheorem{prp}[thm]{Proposition}

\def\End{\mathop{\rm End}\nolimits}

\def\ker{\mathop{\rm ker}\nolimits}

\def\Pic{\mathop{\rm Pic}\nolimits}

\def\deg{\mathop{\rm deg}\nolimits}
\def\rk{\mathop{\rm rk}\nolimits}
\def\tr{\mathop{\rm tr}\nolimits}
\def\td{\mathop{\rm td}\nolimits}

\def\ch{\mathrm {ch}}

\newcommand{\C}{\mathbf{C}}

\newcommand{\Z}{\mathbf{Z}}

\newcommand{\PP}{{\rm P}}
\newcommand{\OO}{{\mathcal{O}}}

\begin{document}
\title{Spinor-valued Higgs fields  \footnote{Based on a talk given at a workshop on June 6th 2023 in Oxford for the FRG Collaborative Research Program ``Complex Lagrangians, integrable systems and quantization".}
}
\author{Nigel Hitchin\\{Mathematical Institute,
Woodstock Road,
Oxford, OX2 6GG}\\{hitchin@maths.ox.ac.uk}}

\maketitle
\begin{abstract}
We investigate the geometry of holomorphic vector bundles $E$ over a Riemann surface $C$ together with a section  
of $\End E\otimes K^{1/2}$ where $K^{1/2}$ is a spin structure -- a square root of the canonical bundle $K$. These parallel to some extent the various features of usual Higgs bundles, such as spectral curve constructions, but some features are radically different. We make essential use of the mod 2 index to distinguish two families of moduli spaces, and provide examples in low genus. 
\end{abstract}

\section{Introduction}\label{intro}
It is now over 35 years   since Higgs fields were introduced in the study of vector bundles on an algebraic curve $C$ \cite{Hitchin1}. Given a vector bundle $E$ a  Higgs field (we can now call this  a ``classical" Higgs field) is a holomorphic section of $\End E\otimes K$. In this article  we replace   $K$ by a choice of square root $K^{1/2}$ -- a spin structure or, in the standard wording of algebraic geometry, a theta characteristic. 

Why? On the one hand a recent physics-inspired paper \cite{Gukov}  introduces Higgs fields twisted with $K^{R/2}$ which suggests further investigation (though perhaps the authors had in mind $R>1$, which makes their calculations for the equivariant Verlinde formula more uniform), on the other it was part of the thesis \cite{Oxbury} of my student Bill Oxbury, written at the same time as  reference \cite{Hitchin1}, and  which I have had occasion to refer to recently. He looked in detail at the case of genus 2 which reveals some of the essential properties. But perhaps the best reason is that it gives a new window onto a familiar subject -- the moduli space $M$ of stable bundles on $C$ -- and focuses on some new features of it and its subspaces.

A key input is parity -- or the  application of the mod 2 index theorem of Atiyah and Singer. Although the Riemann-Roch theorem gives no information about the dimension of $H^0(C,\End E\otimes K^{1/2})$ we do know whether it is odd or even, and this depends on the choice of spin structure. When we consider the situation of $E$ being stable, this introduces two moduli spaces. For the odd case, the generic dimension is 1 and it turns out that we need to study the line bundle $K_M^{1/2}$ over $M$. In the even case the  generic  dimension is zero and we only get a 2-dimensional space by restricting to a hypersurface in $M$. In that case we have instead a rank 2 bundle $V$ over the divisor of a  section of $K_M^{-1/2}$. These two spaces can be regarded as analogues of the cotangent bundle of $M$ which appears in the study of classical Higgs bundles. By using a joint stability condition these spaces extend to larger moduli spaces which have, by the consideration of spectral curves, proper maps to a vector space with generic fibre an  abelian variety. 

An essential difference from the classical case is the appearance of singularities where the dimension of $H^0(C,\End E\otimes K^{1/2})$ jumps by an even number, even when $E$ is stable. We discuss some aspects of the jumping locus but clearly any further investigations about the global moduli space of stable $K^{1/2}$-twisted Higgs bundles will have to  take into  account these singularities. 

We begin in the next section by recalling the  geometry of square roots $K^{1/2}$, then consider the odd and even cases of Higgs fields for stable bundles. In Section \ref{genus2} we look at examples in genus 2 and in the final section make some observations on the jumping locus.  

\section{Spin structures}
Let $C$ be a smooth Riemann surface of genus $g$. A spin structure is defined by a line bundle whose square is isomorphic to the canonical bundle $K$. We denote such a bundle by $K^{1/2}$. If $U$ is a line bundle such that $U^2$ is trivial then the tensor product $UK^{1/2}$ is another square root, and  spin structures form an affine $\Z_2$-space of dimension ${2g}$ with group of translations $H^1(C,\Z_2)$. 

A spin structure in arbitrary dimension describes a vector bundle $V$ and the Dirac operator $\slashed D$ defined on sections of $V$. In two dimensions it is (up to a factor $\sqrt{2}$)
the Cauchy-Riemann operator $\bar\partial: C^{\infty}(K^{1/2})\rightarrow  C^{\infty}(K^{1/2}\bar K)$ and so solutions to the Dirac equation $\slashed D\psi=0$ consist of the space $H^0(C,K^{1/2})$ of holomorphic sections of $K^{1/2}$. 

Given $\psi_1,\psi_2$ smooth sections of $K^{1/2}$, $\psi_1\bar\partial\psi_2$ is a section of $K\bar K$, the bundle of 2-forms, and integrating and applying Stokes' theorem we obtain  
\[\int_C\psi_1\bar\partial\psi_2= -\int_C\psi_2\bar\partial\psi_1\]
so that $\bar\partial$ is formally skew adjoint. Moreover,  a metric in the conformal class is a non-vanishing section $h$ of $K\bar K$ so that multiplication by $h^{-1/2}$ identifies $K^{1/2}\bar K$ with $\bar K^{1/2}$ and then $\bar\partial$ takes sections of $K^{1/2}$ to sections of its conjugate.  The operator is then skew-adjoint antilinear.

The Dirac operator in dimensions  $8k+2$ has this property and then the Atiyah-Singer mod 2 index theorem \cite{Atiyah2} identifies the mod 2 dimension of the kernel with a $KO$-theory characteristic number. In our case this tells us that the parity of $\dim H^0(C,K^{1/2})$ is well-defined regardless of the holomorphic structure on $C$.
Furthermore, of the $2^{2g}$ spin structures $2^{g-1}(2^g-1)$ are odd, meaning the dimension of $H^0(C,K^{1/2})$ is odd. These results are classically derived from the theory of theta-functions \cite{Dol} but it is convenient for us to place them in this wider context. 

Spinors are somewhat mysterious objects but relate to more familiar ones by taking tensor products -- bilinear expressions of spinors are differential forms, a process  that physicists call ``Fierzing". In our case if $s$ is a holomorphic section of $K^{1/2}$ then  $s^2$ is a section of $K$ with double zeros. If we consider the canonical map $C\rightarrow \PP^{g-1}=\PP(H^0(C,K)^*)$ this means that $s$ defines a  hyperplane tangential at each point on the image. 

Thus, for genus $g=2$ the six ramification points of the double cover $C\rightarrow \PP^1$ are the $2^{g-1}(2^g-1)$ odd spin structures, for genus $3$, a plane quartic has $28=2^{g-1}(2^g-1)$ bitangents as  the odd ones and for a generic genus $4$ curve expressed as the intersection of a cubic and a quadric in $\PP^3$ we have $120=2^{g-1}(2^g-1)$ tritangent planes (these last two numbers are famously associated to  the root systems  $E_7,E_8$). 

On the other hand, rather than a plane quartic, a \emph{hyperelliptic} curve of genus 3,   defined by the equation $y^2=(x-x_1)\dots (x-x_8)$ has holomorphic 1-forms with a double zero defined by 
$(ax+b)^2 dx/y$, so in this case one of the 36 even spin structures has a 2-dimensional space of sections. This jumping in dimension is a serious feature which we shall encounter later.  
\section{Higgs fields}
Let $E$ be a rank $n$ holomorphic bundle over $C$. A ``classical" Higgs field is a holomorphic section $\Phi$ of $\End E\otimes K$. A \emph{spinor-valued} Higgs field is a holomorphic section $\Psi$ of $\End E\otimes K^{1/2}$. It is convenient to remove the scalar component of a Higgs field by restricting to the trace zero part  $\End_0E$ and accompanying this by fixing the top exterior power $\Lambda^nE$. In the classical case the Riemann-Roch theorem gives $\dim H^0(C,\End_0 E\otimes K)\ge (n^2-1)(g-1)$ asserting the existence of many Higgs fields, but it gives no information for spinor-valued Higgs fields. On the other hand $\tr (ab)$ is a non-degenerate symmetric form on $\End_0 E$ and 
\[\int_C\tr(\psi_1\bar\partial\psi_2)= -\int_C\tr(\psi_2\bar\partial\psi_1)\]
so that the $\bar\partial$-operator $\bar\partial: C^{\infty}(\End_0 E\otimes K^{1/2})\rightarrow C^{\infty}(\End_0 E\otimes K^{1/2}\bar K)$ is antilinear skew adjoint and the mod 2 index theorem again applies. 
\begin{prp} \label{odd/even} Let $E$ be a rank $n$ holomorphic bundle over $C$ of degree $d$, then if $n$ is odd $\dim H^0(C,\End_0E\otimes K^{1/2})$ is even and if $n$ is even $\dim H^0(C,\End_0E\otimes K^{1/2})=d+\dim H^0(C,K^{1/2})$ mod 2.
\end{prp}
\begin{prf}
Holomorphic structures on $E$ can be viewed as $\bar\partial$-operators $\bar\partial: C^{\infty}(E)\rightarrow C^{\infty}(E\otimes \bar K)$ on the $C^{\infty}$ bundle $E$. These  form an infinite dimensional affine space modelled on $\Omega^{0,1}(\End E)$ and the deformation invariance of the mod 2 index means it is constant on this space. If $E$ has degree $d$ and rank $n$ then the index is the same as for $E=\OO^{n-1}\oplus L$ for a line bundle $L$ of degree $d$ where $\OO^k$ denotes the trivial rank $k$ bundle. Then 
\[ \End_0 E \cong  \OO^{n-1}(L)\oplus \OO^{n-1}(L^*)\oplus \OO^{(n-1)^2}.\]
Now $\dim H^0(C,L^*K^{1/2})=\dim H^1(C,LK^{1/2})$ by Serre duality and from Riemann-Roch  we have 
$d=\dim H^0(C,LK^{1/2})-\dim H^1(C,LK^{1/2})=\dim H^0(C,LK^{1/2})+\dim H^1(C,LK^{1/2}) $ mod 2. 
Thus, modulo 2,  \[\dim H^0(C,\End_0E\otimes K^{1/2})=(n-1)d+ (n-1)^2\dim H^0(C,K^{1/2})\]
and the result follows.
\end{prf}
Given a holomorphic vector bundle $E$ the only spinor-valued Higgs fields may be zero, so we need a means of constructing non-zero examples. We can borrow the spectral curve construction \cite{BNR} from the classical case of $K$-valued Higgs fields. Recall that, in that case, $\det (x-\Phi)=x^n+c_1x^{n-1}+\dots +c_n$ defines a spectral curve $S$  in the total space of $\pi:K\rightarrow C$ and $x$ is the tautological section of $\pi^*K$. Conversely, if $S$ is smooth and $L$ is a line bundle on $S$, its direct image $\pi_*L$ defines a rank $n$ vector bundle $E$ defined by $H^0(U, E)=H^0(\pi^{-1}(U),L)$ for every open set $U\subset C$. Furthermore $x:H^0(\pi^{-1}(U),L)\rightarrow H^0(\pi^{-1}(U),L\pi^*K)$ defines a Higgs field. Clearly we can replace $K$ by $K^{1/2}$ to get a construction of pairs $(E,\Psi)$ where $\Psi$ is a spinor-valued Higgs field. As in the classical case, if we want to fix $\Lambda^nE$ then $L$ must lie not in the full Jacobian of $S$ but in a translate of the Prym variety. 

Now the coefficient $c_i$ is a holomorphic section of $K^{i/2}$ and if $\tr \Psi=0$, $i>1$ and then $\dim H^0(C, K^{i/2})= (i-1)(g-1)$ for $i>2$ and of course $\dim H^0(C, K)=g$. 
 The family of curves then has dimension $g+2(g-1)+\dots +(n-1)(g-1)=1+(g-1)n(n-1)/2$. 
 
 To compute the genus $g_S$ of $S$ note that if $X$ is the total space of $K^{1/2}$ then $K_X\cong \pi^*K^{-1/2}\pi^*K=\pi^*K^{1/2}$. The expression $x^n+c_1x^{n-1}+\dots +c_n$ is a section of $\pi^*K^{n/2}$ on $X$ and so $K_S\cong K_X\pi^*K^{n/2}\cong \pi^*K^{(n+1)/2}$. It follows that $2g_S-2= \deg K_S=n(n+1)(g-1)$ and so the dimension of the Prym variety is given by $g_S-g=n(n+1)(g-1)/2+1-g=(n(n+1)/2-1)(g-1)$. 
 
 Restricting to the smooth spectral curves we have a construction of pairs $(E,\Psi)$ consisting of a family of abelian varieties of dimension $(n(n+1)/2-1)(g-1)$ over an open set in a vector space of dimension $1+(g-1)n(n-1)/2$ and hence a family of dimension $(n(n+1)/2-1)(g-1)+1+(g-1)n(n-1)/2=1+(n^2-1)(g-1).$

 \section{Stable bundles}
 A vector bundle $E$ is stable if for any subbundle $F\subset E$, $\deg F/\rk F<\deg E/\rk E$ and there is a good moduli space $M$ for these which is a projective variety. The only automorphisms of a stable bundle are scalars and so $H^0(C, \End_0 E)=0$ and by Riemann-Roch $\dim H^1(C,\End_0 E)=(n^2-1)(g-1)$ and this space is identified with the tangent space of $M$ at the equivalence class $[E]$ of $E$. The cotangent space is then by Serre duality $H^0(C,\End_0E\otimes K)$.  A generic line bundle on the spectral curve $S$ yields a stable bundle $E$ but it is the pair $(E,\Psi)$ which is constructed so we don't necessarily  see $E$ together with all its Higgs fields.  Note, however, that the family constructed in the previous section has dimension $\dim M+1$.

 We remarked earlier that bilinear expressions in spinors give forms and in the current context there is a natural one: for $\Psi_1,\Psi_2\in H^0(C,\End_0 E\otimes K^{1/2})$ we form the bracket $[\Psi_1,\Psi_2]\in H^0(C,\End_0 E\otimes K)$ and if $E$ is stable this is a cotangent vector to $M$. This has an interpretation:
 \begin{prp} \label{bracket}Suppose $E$ is a stable bundle and  $\alpha\in \dim H^1(C,\End_0 E)$  a tangent direction of the moduli space at $[E]$ along which $\Psi_1,\Psi_2\in H^0(C,\End_0 E\otimes K^{1/2})$ can be extended to first order. Then the cotangent vector $[\Psi_1,\Psi_2]$ annihilates $\alpha$.
  \end{prp} 
 \begin{prf} Let $\dot A\in \Omega^{0,1}(\End_0E)$ be a Dolbeault representative of the class $\alpha$ then a first order deformation of $\Psi_i$ is given by $\dot\Psi_i$ satisfying $\bar\partial \dot\Psi_i+[\dot A,\Psi_i]=0$. Then, using trace for the Serre duality pairing,
 \[  [\Psi_1,\Psi_2](\alpha)= \int_C\tr(\dot A[\Psi_1,\Psi_2])=-\int_C  \tr(\Psi_2[\dot A,\Psi_1]) + \tr(\Psi_1[\Psi_2,\dot A]).  \]
 Rewriting the right hand side using the equation for $\dot\Psi_i$ gives 
 \[\int_C\tr(\Psi_2\bar\partial \dot\Psi_1)-\tr(\Psi_1\bar\partial \dot\Psi_2)= \int_C\bar\partial\tr(\Psi_2\dot\Psi_1-\Psi_1\dot\Psi_2) =0\]
 by Stokes' theorem.
 \end{prf}

 At this point, rather than attempt to construct  a moduli space of stable spinor-valued Higgs bundles -- pairs $(E,\Psi)$ satisfying the stability criterion for $\Psi$-invariant subbundles paralleling the $K$-twisted version \cite{Hitchin1} -- we shall restrict ourselves to considering the situation where $E$ itself is stable, but still bearing in mind the spectral curve construction. There are two cases, where $\dim H^0(C, \End_0E\otimes K^{1/2})$ is odd or even.

  \subsection{The odd case}
When $\dim H^0(C,\End_0 E\otimes K^{1/2})$ is odd we expect the generic value to be one and this provides a line bundle $L$ over an open set in the stable bundle moduli space, realizing 
 $\dim M+1$ from the spectral curve construction. 
 
 This in fact holds, for if, in a neighbourhood of a point  $[E]$ in $M$, the bundles admit two linearly independent Higgs fields then $\Psi_1,\Psi_2$ extend in all directions and so from the proposition $[\Psi_1,\Psi_2](\alpha)=0$ for all $\alpha$. This means $[\Psi_1,\Psi_2]=0$.
  But for a smooth spectral curve $S$,  the generic case,  the constructed Higgs field $\Psi$ is regular at each point and its bundle of centralizers in $\End_0 E$ is isomorphic to the bundle $K^{-1/2}\oplus K^{-1}\oplus \cdots \oplus K^{-(n-1)/2}$ corresponding to the powers $\Psi,\Psi^2,\dots$ So taking $\Psi=\Psi_1$, if $[\Psi_1,\Psi_2]=0$ then $\Psi_2$ is a holomorphic section of $(K^{-1/2}\oplus K^{-1}\oplus \cdots \oplus K^{-(n-1)/2})\otimes K^{1/2}$ and since all terms except the first have negative degree, $\Psi_2$ can only be a multiple of $\Psi_1$. 
  
  The expected codimension of the locus in $M$ where the dimension jumps from 1 to 3 is 3 (indeed  Proposition \ref{bracket} tells us that the three cotangent vectors $[\Psi_1,\Psi_2],[\Psi_2,\Psi_3],[\Psi_3,\Psi_1]$ are conormal at the smooth points of the locus) and the line bundle $L$ extends uniquely to the whole of $M$. In fact we have 
  \begin{prp} \label{root}The line bundle of spinor-valued Higgs fields over the moduli space $M$ of stable bundles is isomorphic to $K_M^{1/2}$.
   \end{prp}   
   \begin{prf}
  There is a universal bundle $\End_0U$ over $M\times C$ with projections $p_1,p_2$ onto the two factors.  Then $R^0_{p_1}(\End_0U\otimes p_2^*K^{1/2})$ is the line bundle $L$ and by Serre duality $R^1_{p_1}(\End_0U\otimes p_2^*K^{1/2})\cong L^*$.  The Grothendieck-Riemann-Roch theorem gives  $\ch(L)-\ch(L^*)=p_{1\ast}(\ch(\End_0U)p_2^*(\ch( K^{1/2})\td(C)))$ but for the second factor   $\ch( K^{1/2})\td(C)=(1+(g-1)x)(1+(1-g)x)=1$ ($x^2=0$ since it is the top degree cohomology class of $C$) so, on the open set where the dimension is one, 
$\ch(L)-\ch(L^*)=p_{1\ast}(\ch(\End_0U))$
 and hence $2c_1(L)=p_{1\ast}(\ch_2(\End_0U))$.
  
  The cotangent bundle $T^*M$ is similarly given by $R^0_{p_1}(\End_0U\otimes p_2^*K)$ and in this case $H^1(C,\End_0U\otimes K)=H^0(C,\End_0U)^*$ which vanishes for a stable bundle so that $R^1_{p_1}(\End_0U\otimes p_2^*K)=0$. Then
  \[ \ch(T^*)=p_{1\ast}(\ch(\End_0U)p_2^*(\ch( K)\td(C)))=p_{1\ast}(\ch(\End_0U)p_2^*(1+(g-1)x))\]
  and so $c_1(T^*)=p_{1\ast}(\ch_2(\End_0U))$ since $\ch_1(\End_0U)=0$. 
  
  It follows from these two calculations that  $c_1(L)=c_1(T^*)/2$.
  But $H^2(M,\Z)\cong \Z$ and $M$ is simply connected so we have an isomorphism of holomorphic bundles $L\cong K_M^{1/2}$.
  \end{prf} 
  If we construct a  moduli space ${\mathcal S}$ of spinor-valued pairs $(E,\Psi)$ using stability defined by $\Psi$-invariant subbundles then we see a parallel here with the classical Higgs bundle moduli space ${\mathcal M}$ which contains $T^*M$ as a dense open subspace. Here we have the similar $K_M^{1/2}\subset {\mathcal S}$.
    
  One  difference concerns the coefficients $c_i$ in the characteristic polynomial. In the classical case a linear form on $H^0(C, K^i)$ defines  a symmetric form on each cotangent 
  space so we obtain $(2i-1)(g-1)$ holomorphic sections of $S^iT$, the $i$th symmetric power of the tangent bundle,  and moreover these generate all of them. By contrast, in the spinor-valued case we get a linear form on $H^0(C, K^{i/2})$ yielding $(i-1)(g-1)$ sections of $K_M^{-i}$, but this is far fewer than given by the Verlinde formula \cite{Thad}. Of course, putting together all choices of spin structure on $C$ may generate more.
  
  Another distinction is the fact that  the line bundle is identified in Proposition \ref{root} by relying on Grothendieck-Riemann-Roch, and not intrinsically -- Serre duality naturally identifies the space of \emph{classical} Higgs fields with the cotangent space of $M$.

\subsection{The even case}\label{even} 
  When $\dim H^0(C,\End_0 E\otimes K^{1/2})$ is even the generic value on the moduli space of stable bundles is zero. This follows  from the argument above: if it is two or more and extends in all directions then $[\Psi_1,\Psi_2]=0$ which is a contradiction. 
  
  Considering the $\bar\partial$-operator on a family of vector bundles $F$ in general, one defines the  Quillen determinant bundle $(\Lambda^{top}H^0(C,F))^*(\Lambda^{top}H^1(C,F))$ and if the index $\dim H^0(C,F)-\dim H^1(C,F)=0$ there is a canonical determinant section which is either zero or vanishes on a divisor where $H^0(C,F)\ne 0$. In our case of $F=\End_0E\otimes K^{1/2}$ where $\bar\partial$ is skew-adjoint, $H^1(C,F)\cong H^0(C,F)^*$ and the determinant is the square of a section of a line bundle, the Pfaffian section. This is the analogue of the fact that the determinant of a skew-symmetric matrix $A$ in even dimensions is the square of a polynomial in $A$. 
     
   But the calculation in the proof of Proposition \ref{root} shows that the determinant bundle is $L^{-2}=K_M^{-1}$ so the Pfaffian is a section of $K_M^{-1/2}$ and vanishes on a hypersurface $P$. For $[E]\in P$ we have generically a 2-dimensional space of sections of $\End_0E\otimes K^{1/2}$ and hence a rank 2 vector bundle $V$. If $\Psi_1,\Psi_2\in H^0(C,\End_0E\otimes K^{1/2})$ we have  $\Psi_1\wedge \Psi_2$ spanning  $\Lambda^2V_{[E]}$ and then $[\Psi_1,\Psi_2]$ gives from Proposition \ref{bracket} an isomorphism from $\Lambda^2V_{[E]}$ to the conormal space. But $P$ is the zero locus of a section of the Pfaffian line bundle, hence  we have $\Lambda^2V\cong K_M^{1/2}$ on $P$.
   
 \begin{rmk}
  Another viewpoint on the vector bundle $V$ is provided by taking a class  $\alpha\in H^1(C,\OO)$ and the natural map
    \begin{equation}
  \alpha: H^0(C, \End_0 E\otimes K^{1/2})\rightarrow H^1(C, \End_0 E\otimes K^{1/2}). \label{disc}
  \end{equation}
  defining a homomorphism from $V$ to $V^*$. As $V$ has rank 2, we can also write $V^*=V\otimes \Lambda^2V^*$
   and given that $\Lambda^2V\cong K_M^{1/2}$, and $\alpha$ is symmetric, this   may be seen as a section of $\End_0V\otimes K_M^{-1/2}$
on $P$, or since $K_PK_M^{1/2}\cong K_M$, a $K_P^*$-twisted Higgs field $\phi\in H^0(P, \End_0V\otimes K^{*}_P)$. 

It follows that  the vector bundle $V$ can be obtained (under generic conditions) as the direct image of a line bundle $U$ on the spectral cover defined by $z^2-\det\phi=0$. This is a hypersurface $X$ in the total space of $\pi:K_P^{*}\rightarrow P$, a double cover of $P$ ramified over a divisor of  $\pi^*(K^{*}_P)$ and so $K_X\pi^*(K_P^{*})\cong \pi^*(K_P^{*})$ and $X$ is a Calabi-Yau variety (which is almost always singular).

Let $\tau$ be the involution on $X$ and $p:X\rightarrow P$ the projection. Then just as in the curve case, the eigenvectors of $p^*\phi$  are interchanged by $\tau$ and coincide where $\det \phi=0$. Then  $U\tau^*U\cong p^*K_M^{-1/2}$. 
\end{rmk}
\subsection{Links with classical Higgs fields}\label{link} 
Here are two cases where spinor-valued Higgs fields link up with the classical case.
\begin{enumerate}
\item If $E$ has rank $n=2m$ suppose the spectral curve $S$ is invariant by the involution $x\mapsto-x$ then it has equation $x^{2m}+a_2x^{2m-2}+\dots + a_{2m}=0$ where $a_{2i}\in H^0(C,K^m)$. Then the quotient, obtained by setting $y=x^2$, lies in the total space of $K$ and becomes a spectral curve $\bar S$ for a classical Higgs field defined by $y^{m}+a_2y^{m-1}+\dots + a_{2m}=0$. Let $L$ be the line bundle on  $\bar S$ for this Higgs bundle then its pull-back is invariant by the involution and the action at the fixed points is the identity. This is an example of the so-called Cayley correspondence but with $K$ replaced by $K^{1/2}$. For $K$-twisted Higgs fields \cite{Schaposnik} \cite{Bradlow}, this defines the spectral curve for a Higgs bundle which corresponds to a flat $U(n,n)$-connection.
\item Take a fixed  $s\in H^0(C,K^{1/2})$ then taking the product with $s$ transforms $H^0(C,\End_0E\otimes K^{1/2})$ to $H^0(C,\End_0E\otimes K)$ so we get an embedding of the spinor-valued moduli space in the classical one. All such Higgs fields vanish at the zeros of $s$ and hence this image lies in the critical locus of the integrable system \cite{Hitchin2}. These are the points in the moduli space of Higgs  bundles where the derivative of the map $h:{\mathcal M}\rightarrow {\mathcal B}$ is not of maximal rank. 
\end{enumerate}

\section{Examples in genus $2$, rank $2$}\label{genus2}
\subsection{Spectral curves}
We consider the genus 2 curve $C$ defined by $y^2=(x-x_1)(x-x_2)\dots(x-x_6)$ and denote by $a_i\in C$ the fixed point of the hyperelliptic involution $\sigma$ lying over $x_i\in \PP^1$. So the odd spin structures are $K^{1/2}\cong \OO(a_i)$. The spectral curve $S$ for a spinor-valued Higgs field is given by $z^2-c=0$ in the total space of $K^{1/2}$ where $c\in H^0(C,K)$. 

The hyperelliptic involution on $C$  acts in different ways on the 2-dimensional space $H^0(C,K^{3/2})$, depending on whether the spin structure is odd or even. In the odd case, if $s\in H^0(C,K^{1/2})$ vanishes at  $a_1$  then each section of $K^{3/2}$ is of the form $sa$ where $a$ is a section of $K$ and the action is by a scalar $\pm 1$. In the even case, $H^0(S,K^{3/2})$ is spanned by sections $s_1,s_2$ with divisors $a_1+a_2+a_3$ and $a_4+a_5+a_6$ respectively. Note that there are ${6\choose 3}/2=10=2^{g-1}(2^g+1)$ ways of choosing these pairs of triples. Here the $\pm $ eigenspaces are spanned by these two distinguished sections. Since $H^0(C,K^{3/2})$ can be identified with the cotangent space of the Prym variety the two actions give different   structures to the Prym varieties.

For the odd case the situation is described in the article \cite{Hitchin2} and the Prym variety corresponding to the differential  $c=(x-a)dx/y$ is isomorphic to the Jacobian of the genus 2 curve $y^2=(x-a)(x-x_2)\dots(x-x_6)$.

 In the even case the Prym variety is a product of elliptic curves: in the equation $z^2-c=0$ for $S$, $z$ is a section of $\pi^*(K^{1/2})$ and   $\pi^*s_1$ is a section of $\pi^*(K^{3/2})$. Then $y_1=zs_1$ satisfies $y_1^2=(x-a)(x-x_1)(x-x_2)(x-x_3)$ defining a map $S\rightarrow E_1$ to an elliptic curve. Similarly $y_2=zs_2$ maps to $E_2$ given by $y_2^2=(x-a)(x-x_4)(x-x_5)(x-x_6)$. Pulling back line bundles from the product of the two curves gives the 2-dimensional Prym variety.

\subsection{$\Lambda^2E$ odd degree}\label{odddegree}
The moduli space $M$ of stable bundles of fixed  determinant $\Lambda^2E$ of odd degree is isomorphic to the intersection of two quadrics in $\PP^5$ \cite{Newstead}. 

\subsubsection{Odd $K^{1/2}$}\label{odd/odd}
Here from Proposition \ref{odd/even} the generic dimension of $H^0(C,\End_0E\otimes K^{1/2})$ is zero and we are concerned with a vector  bundle $V$ over the Pfaffian divisor $P\subset M$. For the intersection of two quadrics, $K_M^{1/2}=\OO(-1)$ for the embedding $M\subset \PP^5$ so $P$, a divisor of $K_M^{-1/2}$, is a hyperplane section of $M$ and hence the del Pezzo surface given by the intersection of two quadrics in $\PP^4$. 
\begin{prp} The vector bundle $V\rightarrow P$ is isomorphic to the cotangent bundle $T^*P$.
\end{prp}
\begin{prf}

 Since  $K^{1/2}$ is odd it has a section $s$ vanishing at $a_1$, say, and  following the approach in Section \ref{link}, $\Psi\mapsto  s\Psi\in H^0(C,\End_0E\otimes K)$ embeds $V$ in $T^*M$. Its  image is an open set in the critical locus\cite{Hitchin2} and in general this is  a symplectic submanifold.

However, we can be more explicit, making use here of the results of a recent paper \cite{BV} which provides a concrete expression for symmetric tensors, sections of  $S^2T$, on the intersection of quadrics in any dimension. These, remarkably, define integrable systems on the cotangent bundle, generalizing the Higgs bundle case of three dimensions. 

If the quadrics are $\sum_{i=1}^n z_i^2=0, \sum_{i=1}^n \mu_iz_i^2=0$ then  the space of symmetric tensors is spanned by the $n$ cases:
 \begin{equation}
 s_i=\sum_{j\ne i}\frac{(z_i\partial_j-z_j\partial_i)^2}{\mu_i-\mu_j} \label{BVeqn}
 \end{equation} 
Here $z_i\partial_j-z_j\partial_i$ is a vector field on the first, standard,  quadric, not tangential to the second but the particular quadratic combination in the formula is well defined as a symmetric tensor on the intersection $M$ of the two.

The embedding $\Psi\mapsto s\Psi$ of $V$ in $T^*M$ gives a 2-dimensional space of classical Higgs fields which vanish at $a_1$. Then the quadratic form $\tr(\Phi_1\Phi_2)(x_1)$ defines a  symmetric form with a 2-dimensional degeneracy subspace. But setting $z_1=0$ in equation \ref{BVeqn} for $i=1$ we obtain
\[s_1=\sum_{j\ne 1}\frac{z_j^2\partial_1^2}{\mu_1-\mu_j} \]
which is clearly of rank $1$. 
Moreover the interior product of this symmetric tensor with a cotangent vector to $z_1=0$ is zero, so the degeneracy subspace is isomorphic to the cotangent space  of  
the divisor of the section $z_1$ of $\OO(1)\cong K_M^{-1/2}$, the Pfaffian locus $P$. The surface $P$  is therefore the 2-dimensional  intersection of quadrics $\sum_{i=2}^n z_i^2=0, \sum_{i=2}^n \mu_iz_i^2=0$  and $V$ is its cotangent bundle. 
\end{prf}
\subsubsection{Even $K^{1/2}$}\label{odd/even}

Suppose now the rank 2 bundle has odd degree, then from Proposition \ref{odd/even}, we have a non-zero element in $H^0(C,\End_0V\otimes K^{1/2})$ if $K^{1/2}$ is one of the $10$ even spin structures. Generically we have a moduli space which is the total space of $K^{1/2}_M$, and in the case of the intersection of two quadrics this is $\OO(-1)$. There is a description \cite{Hitchin3} of this one-dimensional space of Higgs fields  in terms of the quotient of $M$ by the action of $H^1(C,\Z^2)$ by $E\mapsto E\otimes U$ for a line bundle $U$ with $U^2$ trivial. This transformation leaves any Higgs field in $H^0(C,\End_0E\otimes K^{1/2})$ unchanged and so is a natural setting for their study. On the other hand  the quotient space has singularities where $E\cong E\otimes U$. We begin with a description of this quotient.

The intersection of two quadrics $\sum_{i=1}^6 z_i^2=0, \sum_{i=1}^6 \mu_iz_i^2=0$ has an obvious action of $\Z_2^5$ by $z_i\mapsto \pm z_i$ and \cite{Newstead} this is the action of $H^1(C,\Z_2)\cong \Z_2^4$ together with the hyperelliptic involution $\sigma$. Setting $w_i=z_i^2$ gives a map to $\PP^3$ but separating $E$ from $\sigma^*E$ is a double covering branched over the six planes $x_i=0$ and this is the quotient space. 

This description was initially found by Atiyah \cite{Atiyah1} from a different point of view. Each vector bundle has in general 4 maximal subbundles and a choice of one means that the endomorphism bundle  is the same for $E$ given as an  extension $\OO\rightarrow E\rightarrow L$ where $L$ has degree 1. The extension class in the 2-dimensional space $H^1(C, L^*)$ is determined by the section $s\in H^0(C, LK)$ which it annihilates by Serre duality, and the divisor of $s$ is $p+q+r$. The other three subbundles give $(p, \sigma(q),\sigma(r)),(\sigma(p),q,\sigma(r)), (\sigma(p),\sigma(q),r)$. Then the image of $p,q,r\in C$ under the double cover $C\rightarrow \PP^1$ gives a map from the quotient of $M$  to the symmetric product $S^3\PP^1$ which is $\PP^3$. Put differently, this is the projective space of the vector space of polynomials $p(x)$ of degree $3$ and the double covering is branched over the six hyperplanes $p(x_i)=0$.

To construct the Higgs field we map the curve $C$ to $\PP^1\times \PP^1$ using the 2-dimensional space of sections of $K$ on the first factor and of $K^{3/2}$ on the second. As we noted, for an even spin structure the hyperelliptic involution acts non-trivially on the second factor.  This embeds $C$ in $\PP^1\times \PP^1$ which can be identified as a quadric in $\PP^3$, and gives the degree 5 embedding of $C\rightarrow \PP^3$ by sections of $K^{5/2}$.  Then the statement is: 

\begin{prp} \label{pqr} Represent the bundle $E$ as above as an extension $\OO\rightarrow E\stackrel{\pi}\rightarrow L$  defined by  the divisor $p+q+r$ where $p,q,r\in C\subset  \PP^3$ are not collinear. The plane through $\sigma (p),\sigma (q),\sigma (r)$ meets $C$ in two further points, the divisor of a section $c$ of $LK^{1/2}$. Then there is a unique Higgs field  $\Psi\in H^0(C,\End_0E\otimes K^{1/2})$ such that $\pi\Psi(1)=c$.
\end{prp} 
The proof \cite{Hitchin3} is a simple Dolbeault argument.
\begin{rmk}\label{collinear}
What if the points are collinear? Then $L K^{1/2}$, which  has degree $2$ and a two-dimensional space of sections  must be isomorphic to $K$, or equivalently $L\cong K^{1/2}$ and $p+q+r$ is a divisor of $K^{3/2}$. In this case the extension is $\OO\rightarrow E\rightarrow K^{1/2}$ and there is a nilpotent Higgs field $e\mapsto \pi(e)\in \ker \pi \otimes K^{1/2}$. Conversely every  bundle with a nilpotent Higgs field is an extension of this form, describing a   projective line $\PP(H^1(C,K^{-1/2}))$ in the quotient space  of $M$ and 16 lines in $M$ itself.
\end{rmk}

\subsection{ $\Lambda^2E$ trivial}
The moduli space $M$ of (semi)-stable bundles is in this case projective space $\PP^3$ \cite{NR}. To each  bundle $E$  the  line bundles $L$ of degree $1$ such that $H^0(C,E\otimes L)\ne 0$ describe a curve in $\Pic^1(C)$ which is the divisor of a section in the linear system $2\Theta$  and hence a point in $\PP(H^0(\Pic^1(C), 2\Theta))=\PP^3$. Strictly semistable bundles are represented by their S-equivalence class,   a direct sum  $U\oplus U^*$ where $\deg U=0$. In the moduli space $\PP^3$  this is a Kummer quartic surface, a quotient of $\Pic^0(C)$. 

If $U$ is a line bundle with $U^2$ trivial  then $E\mapsto E\otimes U$  again gives an action of $H^1(C,\Z_2)$ on $\PP^3$. A central extension by $\Z_2$ acts on the 4-dimensional vector space $H^0(\Pic^1(C), 2\Theta)$ and in a standard theta-function basis for these sections the action is generated by double transpositions and double sign changes. 
\subsubsection{Odd $K^{1/2}$}
 Take  the spin structure $K^{1/2}=\OO(a_1)$,  then  $H^0(C,\End_0E\otimes K^{1/2}) \ne 0$ from Proposition \ref{odd/even} and we have pairs $(E,\Psi)$ defined by points on the total space of the line bundle $K_M^{1/2}$ on $M=\PP^3$. Since $K_M=\OO(-4)$,  $K_M^{1/2}=\OO(-2)$ and $\Psi\mapsto s\Psi$ embeds the total space of $K_M^{1/2}$ as 
a distinguished subbundle $\OO(-2)\subset T^*\PP^3$. 

We can identify it by using an explicit description of the integrable system \cite{Gaw}, and the procedure in the paper \cite{Hitchin3}. The cotangent bundle of $\PP^3$ is expressed  as a complex symplectic quotient  
$\{(p,q)\in \C^4\times \C^4: q\ne 0, \sum_{i=1}^4p_iq_i=0\}/\C^*$ where the action is $\lambda.(p,q)=(\lambda^{-1}p,\lambda q)$. So $(q_1,\dots,q_4)$ are homogeneous coordinates on $\PP^3$. Then the map $h:T^*\PP^3\rightarrow H^0(C,K^2)$ is
\[\sum_{i\ne j}\frac{r_{ij}(q,p)}{(x-x_i)(x-x_j)}dx^2=\sum_{i\ne j}\prod_{k\ne i,j}(x-x_k){r_{ij}(q,p)}\frac{dx^2}{y^2}\]
where $r_{12}(q,p)= (q_1p_1 + q_2p_2-q_3p_3-q_4p_4)^2$ and similar terms. 

Then, as in Section \ref{odd/odd}, when $\Phi_0=\Psi s$, $\tr(\Phi\Phi_0)$ vanishes at $x_1$ for all $\Phi$, so the quadratic form $\tr(\Phi_1\Phi_2)(x_1)$ is degenerate and $\OO(-2)$ is the degeneracy subspace.
Now evaluating the expression above at $x=x_1$  gives a linear combination of the terms    $r_{1i}$. Each of these is a quadratic form $\ell_i\otimes \ell_i$ where $\ell_i(p_1,p_2,p_3,p_4)$ is a linear form. But on examination one  sees that $\ell_i(q_2, -q_1, q_4, -q_3)=0$ for each $i$, so   the degeneracy subspace, and hence the image of $K_M^{1/2}$ in $T^*M$, is $(p_1,p_2,p_3,p_4)=(q_2, -q_1, q_4, -q_3)$. 

\subsubsection{Even $K^{1/2}$}
From Proposition \ref{odd/even} $\dim H^0(C,\End_0E\otimes K^{1/2})$ is even and  the stable part of the moduli space is a rank $2$ bundle $V$ over a Pfaffian divisor of $K_M^{-1/2}=\OO(2)$, namely a quadric surface $Q$ in $\PP^3$. Here a choice of spectral curve gives a well-defined  degree 4 map from the Prym variety, a product of two elliptic curves, to the   quadric surface, a product $Q\cong \PP^1\times \PP^1$.

In this case, since $K^{1/2}$ is even and has no sections, we cannot embed it into $T^*M$ and so we need other methods to determine the bundle $V$. We know that $\Lambda^2V\cong  K^{1/2}_M$ and hence $c_1(V)=-c_1(TM)/2$. 

Although we can't embed this vector bundle into $T^*M$ as in the odd degree case, nevertheless we have:
\begin{prp} The vector bundle $V\rightarrow Q$ is isomorphic to the cotangent bundle $T^*Q$.
\end{prp} 

\begin{prf}
 The argument is a roundabout one, 
using a recent classification of nef vector bundles on a quadric surface \cite{Ohno}. We first show that $V^*$ is nef, meaning the line bundle $H$ on the projective bundle $\PP(V)$, dual to the tautological line bundle, has non-negative degree on any curve. If $H$ is not nef, then $c_1(H)[D]<0$ for a curve $D$ which is in the base locus of the linear systems $H^n$ for all $n>0$. 

The information we have about $V$ is the quadratic map $V\rightarrow H^0(C,K)$ which defines a 2-dimensional space of sections of $S^2V^*$, or equivalently  sections of  $H^2$ on $\PP(V)$. These are parametrized by the dual of $H^0(C,K)$ and this  is equivalent in genus 2 to evaluation of $\tr\Psi^2$ at a point $x\in \PP^1$. We consider the base locus of this pencil of surfaces $\{S_x\subset \PP(V): x\in \PP^1\}$, namely the  intersection $S_x\cap S_y$. 

The surface $S_x$ 
meets a fibre of $\PP(V)\rightarrow Q$ in two, possibly coincident, points which are the null directions of the  quadratic form $\tr \Psi^2(x)$ and so a point in the intersection of $S_x$ for all $x$ is a nilpotent Higgs field. This expresses $E$ as an extension $L\rightarrow E\rightarrow L^*$, with $L$ the kernel of $\Psi$ and $\Psi$ given by  a section of the degree $1$  line bundle $L^2K^{1/2}$.  So $L^2K^{1/2}\cong \OO(x)$ for some $x\in C$. This is a strictly semistable bundle, and represented in the moduli space by a direct sum $L\oplus L^*$,  and so the locus lies in the intersection of the Kummer surface in $\PP^3$ with $Q$. (Note that  $\End_0E\otimes K^{1/2}\cong L^2K^{1/2}\oplus K^{1/2}\oplus L^{-2}K^{1/2}\cong \OO(x)\oplus K^{1/2}\oplus \OO(\sigma(x))$ has a 2-dimensional space of sections since $K^{1/2}$ is an even spin structure.) 

The bundles $L, L^*$ define the same vector bundle $E$, so the unordered pair $(x,\sigma(x))\in \PP^1=C/\sigma$ defines the Higgs field. Choosing $L,L^*$ from $L^2\cong K^{-1/2}(x)$ gives a 16-fold covering of $\PP^1$ branched over the six points $x_i$, where the inverse image is 8 points, so by Riemann-Hurwitz this gives a smooth genus $9$ curve in the intersection of the Kummer surface with $Q$. But this intersection is a curve in the linear system $(4,4)$ on $Q=\PP^1\times \PP^1$ and the genus of a smooth member is $(4-1)^2=9$ so the intersection is a connected smooth curve. In all cases we have two null vectors so the inverse image in each surface $S_x$ is a curve $D$ which is an unramified double cover and hence has Euler characteristic
$2(2-2\times 9)=-32$.   So $D$ is the base locus of the pencil.

A smooth  divisor $S$ of a section of $H^2$ has $K_S\cong K_{\PP(V)}H^{2}$ and $K_{\PP(V)}\cong p^*(\Lambda^2V^*)H^{-2}p^*(K_Q)$. But $\Lambda^2V\cong K_Q$ so $K_S$ is trivial. Then the self-intersection of $D$ in $S$ is $D^2=2c_1(H)[D]=32$. This is positive hence the bundle $H$ is nef.

Now Ohno \cite{Ohno} shows that, given the first Chern class $c_1=\OO(2,2)$ (our case for $V^*$), the rank $2$ vector bundles on $\PP^1\times \PP^1$ which are nef are:
\[\OO(1,1)\oplus \OO(1,1)\quad \OO(2,2)\oplus \OO\qquad \OO(2,0)\oplus \OO(0,2)      \]
and the last one is the tangent bundle. They are distinguished by the value of the second Chern class $c_2$. But we have the relation $h^2+c_1(V)h+c_2(V)=0$ where $h=c_1(H)$
and so $h^3[\PP(V)]=(c_1(V)^2-c_2(V))[Q]$.  But $(2h)^3[\PP(V)]=D^2=32$. So $4=h^3[\PP(V)]=8-c_2(V)$ and $c_2(V)=4$ which is the cotangent bundle. 
\end{prf}

\section{Jumping loci}
We have already seen the dimension of $H^0(C,\End_0E\otimes K^{1/2})$ jump from $0$ to $2$ on a hypersurface, but this is just the first of the jumps in dimension as we vary in the moduli space of stable bundles, or more generally in any moduli space of stable spinor-valued Higgs bundles. Such jumping occurs for  $\dim H^0(C,\End_0E\otimes K)$ in the classical Higgs bundle moduli space but this is frequently compensated for by the corresponding jump in $\dim H^0(C,\End_0E)$,  tangential to the gauge transformations which give equivalence.

In general, these jumps can be detected topologically, in our case by treating the $\bar\partial$-operator as being part of a family of Fredholm operators.
The space of Fredholm operators on a Hilbert space has different components corresponding to the index and in \cite{Ko} Koschorke identified finite codimensional subspaces defined by the condition $\dim \ker\ge n$ and identified the cohomology class of these.  Pulling back to a finite-dimensional family gives a cohomology class which generically is represented by the subspace of operators for which the null-space has dimension $\ge n$ -- a jumping locus. For a compact family parametrized by $M$ the kernel and cokernel are represented by virtual vector bundles $V,U$ over $M$ and the cohomology class is given by a universal formula in the characteristic classes appearing in $\ch(V)-\ch(U)$. 

Our 2-dimensional Dirac operator is skew adjoint and antilinear or, using the metric, is formally a skew-symmetric Fredholm operator  from  a space to its dual. Koschorke's approach can be adapted to this situation where the finite-dimensional model is of a skew-symmetric homomorphism $A$ from a vector bundle $V$ to its dual $V^*$. There are two cases where $\rk V$ is odd and $\dim \ker A\ge 1$ always and the other case $\rk V$ even. These correspond to the mod 2 index being 1 or 0 respectively. 

In the even case the jump from generic zero to $\ge 1$ is clear: $A:V\rightarrow V^*$ where $V$ has rank $2m$ is a section $\alpha$ of $\Lambda^2V^*$ and is degenerate if $\alpha^m\in \Lambda^{2m}V^*$ vanishes. The cohomology class is $c_1(V^*)$. This is  half the Chern class of the  usual determinant line bundle, and we have encountered this already. In general, the jump is represented \cite{HT}  by classes in $c_i(V^*)$
 \[ c_1 \qquad \begin{vmatrix}
c_2 & c_3\\
1 & c_1 \\
\end{vmatrix}\qquad
 \begin{vmatrix}
c_3 & c_4 & c_5 \\
c_1 & c_2 & c_3 \\
0 & 1 & c_1 \\
\end{vmatrix}\qquad
\begin{vmatrix}
c_4 & c_5 & c_6 & c_7\\
c_2 & c_3 & c_4 & c_5 \\
1 & c_1 & c_2 & c_3 \\
0 & 0 & 1 & c_1 
\end{vmatrix}\qquad
\cdots
\]
in codimension $1,3,6,10,...$ the expected dimension of the linear combinations $[\Psi_i,\Psi_j]\in H^0(C, \End_0E\otimes K)$ which give conormal vectors. 
\begin{rmk}
\begin{enumerate}
\item
Note that although these are polynomials in the Chern classes of $V^*$, the interpretation means  they are also defined in terms of  $\ch(V)-\ch(V^*)$: the odd degree components of $\ch(V^*)$. The first three cases are: $\ch_1, \ch_1^3/3-2\ch_3, \ch_1^6/45 - 2 \ch_1^3 \ch_3/3 - 4 \ch_3^2 + 24 \ch_1 \ch_5$.
\item
 In the case of genus $2$ and $\Lambda^2E$ trivial and an odd spin structure the jump occurs in codimension 2 rather than the expected $3$: when $E=L\oplus L^*$ then 
$\End_0 E\otimes K^{1/2}= L^2K^{1/2}\oplus K^{1/2}\oplus L^{-2}K^{1/2}$ so if $L^2K^{1/2}\cong \OO(x)$ then we get a 3-dimensional space. This is nevertheless on the semistable locus which exceptionally is a smooth part of the moduli space. 
\end{enumerate}
\end{rmk} 

We want to apply these formulas  to the case where $M$ is the moduli space of stable bundles  of rank $2$, odd degree and  fixed determinant -- a smooth compact manifold of complex dimension $3g-3$. We already calculated  
$$\ch(V)-\ch(V^*)=\pi_*(\ch(\End_0U)\ch(K^{1/2})\td(C))=\pi_*(\ch(\End_0U)).$$
The rational cohomology is generated \cite{New} by  classes $\alpha\in H^2(M,\Z), \psi_i\in H^3(M,\Z)$ where $1\le i\le 2g$ and $\beta\in H^4(M,\Z)$ where the second Chern class  of the universal bundle $\End_0U$ over $M\times C$ is expressed in K\"unneth components as
\begin{equation}\label{c2}
c_2(\End_0(U))=  2\alpha f -\beta+4\sum_{i=1}^{2g}\psi_ie_i.
\end{equation}
Here $f$ is the fundamental class of $C$, and $e_i$ form a  symplectic basis for $H^1(C,\Z)$ with $e_ie_{i+g}=-f$. 
We calculate $\ch(\End_0(U))$ from (\ref{c2}) and then 
\begin{equation}\label{sinhb}
\pi_*\ch(\End_0(U))=-2\alpha\frac{\sinh \sqrt{\beta}}{\sqrt{\beta}}+2\gamma\frac{\cosh \sqrt{\beta}}{\beta}-2\gamma\frac{\sinh  \sqrt{\beta}}{\beta ^{3/2}}
\end{equation}
where $\gamma=2\sum_1^g\psi_i\psi_{i+g}$.  Since $\ch_1(T)=2\alpha, \ch_2(T)=(g-1)\beta, 3\ch_3(T)=\alpha\beta-4\gamma$\cite{New}, this can also be expressed in terms of Chern classes of the tangent bundle.

We apply this to the moduli space $M$ of Section \ref{odddegree}, the intersection of two quadrics.
\begin{prp} Let $M$ be the 3-dimensional moduli space of stable bundles $E$ with $\Lambda^2E$ fixed and odd degree on a curve $C$ of genus $2$, then  
\begin{enumerate}
\item
 the jumping locus  for $H^0(C,\End_0E\otimes K^{1/2})$ where $K^{1/2}$ is an even spin structure consists of 8 points,
 \item
 there are two orbits under the action of $H^1(C,\Z_2)$ on $M$ and these are represented in the quotient by the two divisors $a_1+a_2+a_3, a_4+a_5+a_6$ of $K^{3/2}$ where $a_i$ are the fixed points of the hyperelliptic involution.
\end{enumerate}
\end{prp} 
\begin{prf}
 
Represent $\End_0E$ by an extension $\OO\rightarrow E\stackrel{\pi}\rightarrow L$ as in Section \ref{odd/even}, then any Higgs field $\Psi$ defines $c=\pi\Psi(1)\in H^0(C,LK^{1/2})$. But  $LK^{1/2}$ has degree 2 and so at most a 2-dimensional space of  sections. So if  $\dim H^0(C,\End_0E\otimes K^{1/2})=3$ there must exist a nilpotent Higgs field which, as in Remark \ref{collinear}, means that $L\cong K^{1/2}$ and the jumping locus is at most one-dimensional. 

Now  for $x\in C$, $H^0(C,\End_0E\otimes K^{1/2})\cong (\End_0E\otimes K^{1/2})_x$ by restriction since $K^{1/2}(-x)$ has degree zero and by stability $H^0(C,\End_0E\otimes K^{1/2}(-x))=0$. This implies that $[\Psi_1,\Psi_2],[\Psi_2,\Psi_3],[\Psi_3,\Psi_1]$ are linearly independent conormal vectors hence the locus is zero-dimensional. 

 In three dimensions we now use from equation \ref{sinhb} $\pi_*\ch(\End_0(U))=-2\alpha +(-\alpha \beta+4 \gamma)/3$. 
 The jumping locus then defines the cohomology class
\[(c_1c_2-c_3)(V^*)=\frac{1}{3}(\alpha^3-\alpha\beta+4\gamma)=\frac{1}{3}(\ch_1(T)^3/8-3\ch_3(T)).\]
For the intersection of two quadrics $T\oplus\OO(2)\oplus\OO(2)\cong TP^5\vert_M$ so $\ch(T)=6e^h-2e^{2h}$ and $\ch_1(T)=2h, \ch_3(T)=-8h^3/3$ where $h\in H^2(M,\Z)$ is the hyperplane class. Hence $(c_1c_2-c_3)(V^*)=2h^3[M]=8$ since the intersection of two quadrics has degree 4.

Then the action of $H^1(C,\Z_2)$ of order 16 on the 8 points in  $M$ implies the existence for each point of an isomorphism  $\alpha: E\rightarrow E\otimes U$ for some line bundle $U$ with $U^2$ trivial, and then $\pi\alpha(1)$ is a non-zero section of $K^{1/2}U$, vanishing at a divisor  $a_i$ of an odd spin structure.   The extension class in $H^1(C, K^{-1/2})$ then arises in the exact sequence 
$$\cdots \rightarrow H^0(a_i, U) \rightarrow H^1(C, K^{-1/2})\rightarrow H^1(C,U)\rightarrow \cdots $$ 
from a section of $U$ on $a_i$. Both are fixed by the hyperelliptic involution $\sigma$ and so the extension class is fixed by $\sigma$. Since  $H^1(C, K^{-1/2})$ is 2-dimensional the projective spaces $\PP(H^1(C, K^{-1/2})), \PP(H^0(C,K^{3/2}))$ are canonically isomorphic  
and  the non-trivial action of $\sigma$ gives just two fixed points in the projective line $\PP(H^0(C,K^{3/2})$, the sections with divisors $a_1+a_2+a_3, a_4+a_5+a_6$.
\end{prf}

Oxbury\cite{Oxbury} gives a geometric argument for this result, and a description of the 16 inverse images in $M$ of the line $\PP(H^0(C,K^{3/2}))\cong \PP(H^1(C, K^{-1/2}))$ in the quotient. Each 
line has a distinguished pair of points which lie in the jumping locus but the lines meet in fours at these points giving $8=(2\times 16)/4$ as calculated from the formula above.


\begin{thebibliography}{9}
  
  

\bibitem{Atiyah1}
   M.~F. Atiyah, Complex fibre bundles and ruled surfaces,
   \emph{Proc. Lond. Math. Soc.}, {\bf 5}, 407--434 (1955).
   \bibitem{Atiyah2}
    M.~F. Atiyah and I.~M. Singer, The index of elliptic operators V,  \emph{Ann. of Math.}, {\bf 93}, 139--149  (1971).
   \bibitem{BNR}
   A. Beauville, M.~S. Narasimhan and S.Ramanan, Spectral curves and the generalized theta divisor, \emph{J.reine angew. Math.} {\bf 398}, 169--179 (1989).
   \bibitem{BV}
   A. Beauville, A. Etesse, A. H\"oring, J. Liu and  C. Voisin,
{Symmetric tensors on the intersection of two quadrics and Lagrangian fibration},  arXiv:2304.10919.  
\bibitem{Bradlow}
S. Bradlow, O. Garcia-Prada and P. Gothen,
{Surface group representations and $U(p,q)$-Higgs bundles}, \emph{J. Differential Geometry,} {\bf 64}, 111-170 (2003).
\bibitem{Dol}
I.~V. Dolgachev, \emph{Classical algebraic geometry}, Cambridge University Press, Cambridge (2012).
  \bibitem{Gaw}
K. Gaw\c{e}dzki  \& P. Tran-Ngoc-Bich, 
Self-duality of the $SL_2$ Hitchin integrable system at genus $2$,
\emph{Commun.Math.Phys.}, {\bf 196},   641--670 (1998).
   \bibitem{Gukov}
   S. Gukov, A. Sheshmani and S-T. Yau,  {3-manifolds and Vafa--Witten theory}, arXiv:2207.05775. 
   \bibitem{HT}
J.Harris and L.Tu, On symmetric and skew-symmetric determinantal varieties, \emph{Topology},  {\bf 23}, 71--84 (1984). 
    \bibitem{Hitchin1}
    N.~J. Hitchin,  The self-duality equations on a Riemann surface, \emph{Proc. London
Math. Soc.},  {\bf 55},  59--126 (1987).
 \bibitem{Hitchin2}
 N.~J. Hitchin, Critical loci for Higgs bundles, \emph{Commun. Math. Phys.}, {\bf 366}, 841--864 (2019).
  \bibitem{Hitchin3}
 N.~J. Hitchin, A note on coupled Dirac operators, \emph{SIGMA}, {\bf 19}, 003 (2023).
 \bibitem{Ko}
U. Koschorke,  Infinite dimensional K-theory and characteristic classes of Fredholm bundle maps, \emph{ Proc. Symp. Pure Math. Amer. Math. Soc.}, {\bf 15},  95--133 (1970). 
 \bibitem{NR}
 M.~S. Narasimhan and S. Ramanan,
Moduli of vector bundles on a compact Riemann surface,
\emph{Ann. of Math.} {\bf 89},  14--51 (1969). 
 \bibitem{Newstead}
P. Newstead, Stable bundles of rank 2 and odd degree over a curve of genus 2, \emph{Topology}, {\bf 7},  205--215 (1968). 
\bibitem{New}
P. Newstead, {Characteristic classes of stable bundles of rank $2$ over an algebraic curve}, \emph{Trans. Amer. Math.Soc.,} {\bf 169}, 337--345 (1972). 
\bibitem{Ohno}
M. Ohno, Nef vector bundles on a quadric surface with first Chern class $(2,2)$, arXiv: 2311:02830.
    \bibitem{Oxbury}
   W. Oxbury, Stable bundles and branched coverings over Riemann surfaces, \emph{D.Phil thesis}, Oxford (1987).
   \bibitem{Schaposnik}
   L.~P. Schaposnik,  Spectral Data for $U(m, m)$-Higgs Bundles, \emph{International Mathematics Research Notices},  {\bf 11}, 3486--3498 (2015).
 \bibitem{Thad}
   M. Thaddeus,  Stable pairs, linear systems and the Verlinde formula, \emph{Inventiones Math.}, {\bf 117}  317--353  (1994).

\end{thebibliography}
\end{document}